\documentclass{amsart}

\usepackage{amsmath,amssymb}
\usepackage{color}

\newtheorem*{thm}{Theorem}
\newtheorem{lem}{Lemma}
\newtheorem*{prop}{Proposition}
\newtheorem{cor}{Corollary}

\theoremstyle{remark}
\newtheorem{rk}{Remark}

\begin{document}

\title{A Criterion for being a Teichm\"uller curve}

\author[E. Goujard]{Elise Goujard}
\address{IRMAR, Universit\'e de Rennes 1, Campus de Beaulieu, 35042, RENNES, FRANCE}
\email{elise.goujard@univ-rennes1.fr}

\date{}

\dedicatory{}

\begin{abstract}
Teichm\"uller curves play an important role in the study of dynamics in polygonal billiards. 
In this article, we provide a criterion similar to the original M\"oller's criterion, to detect whether a complex curve, embedded in the moduli space of Riemann surfaces and endowed with a line subbundle of the Hodge bundle, is a Teichm\"uller curve, and give a dynamical proof of this criterion.

\end{abstract}

\maketitle
\section{Introduction}

Given  a  curve  in  the moduli space of Riemann surfaces, we want to
know  whether  it  is  a Teichm\"uller curve. By Deligne semisimplicity
theorem  the Hodge bundle over the curve decomposes into a direct sum
of  flat  subbundles  admitting variations of complex polarized Hodge
structures  of  weight  $1$.  Suppose  that  the  restriction  of the
canonical   pseudo-Hermitian  form  to  one  of  the  blocks  of  the
decomposition has rank $(1,r-1)$. We establish an upper bound for the
degree  of  the corresponding holomorphic line bundle in terms of the
(orbifold)  Euler  characteristic  of the curve. Our criterion claims
that  if  the  upper  bound  is  attained, the curve is a Teichm\"uller
curve.

For   those  Teichm\"uller  curves  which  correspond  to  strata  of
\textit{Abelian}   differentials   our  criterion  is  necessary  and
sufficient  in  the sense that if the curve is a Teichm\"uller curve,
then  the  decomposition  of  the Hodge bundle necessarily contains a
nontrivial  block  of  rank $(1,1)$ corresponding to the tautological
line bundle for which the upper bound is attained.

The  original  criterion  in  the  same  spirit  was  found by Martin
M\"oller  in  \cite[Th.  2.13  and  5.3]{moller}  where the condition
detecting a Teichm\"uller curve is  formulated  in terms of Higgs bundle, or equivalently in terms of non-vanishing  of  the second fundamental form (Kodaira-Spencer map).
In \cite{wright}, A.~Wright gives an alternative  version  of M\"oller's criterion, in terms of non-vanishing of the period map.

The  key idea of our criterion is based on Forni's observation
that the    tautological    bundle    on   a   Teichm\"uller   curve is spanned by those vectors of the Hodge bundle which  have  the maximal variation  of the Hodge norm  along  the  Teichm\"uller  flow.

We  combine  this  result of Forni with the Bouw--M\"oller version of
the  Kontsevich  formula for the sum of the Lyapunov exponents of the
Hodge  bundle  along  the Teichm\"uller geodesic flow. Similar to the
criteria  mentioned  above,  the  fact  that the Teichm\"uller metric
coincides with the Kobayashi metric will be crucial for the proof.

\section{Criterion}

Having   a   Riemann   surface   $X$,  the  natural  pseudo-Hermitian
intersection   form   on   $H^1(X,   \mathbb   C)$,   is  defined  on
closed 1-forms representing cohomology classes as:
$$
(\omega_1,   \omega_2)=\frac{i}{2}\int_X \omega_1\wedge\overline{\omega_2}\,.
$$
Restricted to $H^{1,0}(X, \mathbb C)$, the form is positive-definite,
and   restricted   to   $H^{0,1}(X,   \mathbb   C)$,   the   form  is
negative-definite.

This pseudo-Hermitian form of signature $(g,g)$ induces a form on the
Hodge  bundle  $H^1$  over  the  the  moduli  space $\mathcal M_g$ of
Riemann surfaces of genus $g$, where the fiber $H^1_{X}$ of the Hodge
bundle over a point $X$ in $\mathcal M_g$ is $H^1(X, \mathbb C)$. The
pseudo-Hermitian  form  is  covariantly  constant with respect to the
Gauss--Manin flat connexion on the Hodge bundle.

Let  $\mathcal  C$  be  a complex curve in $\mathcal M_g$. We want to
detect,  whether  $\mathcal  C$  is  a  Teichm\"uller  curve  or not.
Throughout  this  paper  we  assume  that  the  genus $g$ is strictly
greater  than $1$ since in genus one the problem becomes trivial: the
moduli space $\mathcal M_{1,1}$ is a complex curve itself.

By Deligne semisimplicity theorem \cite[Prop. 1.13.]{deligne}, the Hodge bundle over $\mathcal C$
splits into a direct sum of orthogonal flat subbundles, such that the
restriction  of the canonical pseudo-Hermitian form to each subbundle
is   nondegenerate.  Assume  that  this  splitting  contains  a  flat
subbundle  $\mathbb  L$  of  rank $r$, where $r\geq 2$, such that the
signature  of  the  canonical  pseudo-Hermitian  form  restricted  to
$\mathbb  L$  is  $(1,r-1)$. Let us define $\mathbb L^{1,0}=\mathbb L
\cap H ^ {1,0}$ and $\mathbb L^{0,1}=\mathbb L \cap H^{0,1}$. Deligne
semisimplicity  theorem combined with our assumption on the signature
implies  that  $\mathbb  L  ^{1,0}$ is a holomorphic line bundle over
$\mathcal C$.

Note  that  for  Teichm\"uller  curves corresponding to the strata of
Abelian  differentials  the  splitting is always nontrivial, since it
contains  a  flat  subbundle of rank $2$ such that the restriction of
the  pseudo-Hermitian  form  to this subbundle has signature $(1,1)$.
The  corresponding line bundle $\mathbb L ^{1,0}$ is the tautological
line bundle over the Teichm\"uller curve.

The  curve $\mathcal C$ may have a finite number of cusps and conical
points,  so  we  need  to  consider  the  Deligne  extension  of  the
holomorphic    line    bundle    $\mathbb    L^{1,0}$,   denoted   by
$\overline{\mathbb L^{1,0}}$: it becomes an orbifold vector bundle at
the  cusps and conical points. So it has an orbifold degree, which in
general is not an integer, but a rational number.

Let  $\chi(\mathcal  C)$  be  the generalized Euler characteristic of
$\mathcal   C$:  it   is   given   by  the  formula  $$\chi(\mathcal
C)=2-2g-n_{\mathcal  C} + \sum_i (k_i-1),$$ where $n_{\mathcal C}$ is
the number of cusps on $\mathcal C$, and $2\pi k_i$ is the cone angle
of  the  $i$-th  conical  point.

\begin{thm}
 If  $\chi(\mathcal C)\ge 0$, then $\mathcal C$ is not a Teichm\"uller
curve.

Suppose  that  $\chi(\mathcal C)< 0$. 
For any flat subbundle $\mathbb
L$  of  the  Hodge  bundle  over  $\mathcal  C$  satisfying the above
assumptions, one has
\begin{equation}\label{arakelov}
\deg\overline{\mathbb L^{1,0}}\le-\frac{\chi(\mathcal   C)}{2}\,.
\end{equation}

If  the  equality  is  attained,  then  $\mathcal C$ is a Teichm\"uller
curve, and the line bundle $\mathbb L^{1,0}$ is the tautological bundle.

 Any  Teichm\"uller  curve  corresponding  to  a  stratum  of  Abelian
differentials admits a flat subbundle $\mathbb L$ of the Hodge bundle
satisfying the above conditions, such that
$$
\deg\overline{\mathbb L^{1,0}}=-\frac{\chi(\mathcal   C)}{2}\,.
$$
\end{thm}

The first statement of the theorem results from the fact that any Teichm\"uller curve has negative curvature, or equivalently, an orbifold genus strictly greater than 1. So from now on, we assume that $\chi(\mathcal C)<0$, that is, $\mathcal C$ is hyperbolic. 

\begin{rk}
We  have  to  admit  that  our  criterion  does  not  directly detect
Teichm\"uller   curves  corresponding  to  the  strata  of  quadratic
differentials.  However,  the  canonical double covering construction
associates   to  every  such  Teichm\"uller  curve  $\mathcal  C$  in
$\mathcal   M_g$   a   Teichm\"uller   curve  $\mathcal{C}'$  in
$\mathcal{M}_{g'}$ in the moduli space of curves of larger genus,
such  that  the  new  Teichm\"uller curve already corresponds to some
stratum  of Abelian differentials, and thus, would be detected by our
criterion.   The   new  Teichm\"uller  curve  $\mathcal{C}'$  is
isomorphic to the initial curve $\mathcal C$.

Since  any  Riemann surface $X'$ in the family $\mathcal C'$
admits  a  holomorphic  involution  which  changes  the  sign  of the
corresponding   Abelian   differential,   the   Teichm\"uller  curves
corresponding   to   such  a  double  covering  construction  can  be
identified.  Thus, indirectly the criterion detects all Teichm\"uller
curves.
\end{rk}

\begin{rk} 
As Martin M\"oller pointed out to the author, in the case $r=2$ this theorem can be refound by algebraic methods. Inequality (\ref{arakelov}) is a specific case of Arakelov inequality (see e.g. \cite[Lemme 3.2]{deligne}), and the bound is attained if and only if $\mathbb L$ is maximal Higgs (see \cite{pe2} or \cite{vz} for generalization in higher dimension). So M\"oller's criterion applies here and gives the conclusion of the theorem.
\end{rk}

\section{Comparison of hyperbolic versus Teichm\"uller metric}

Recall  that  at  any  point  $X\in\mathcal  M_g$  the  tangent space
$T_X\mathcal  M_g$  is identified with the space of essentially bounded
Beltrami  differentials,  which  is  in  duality  with  the  space of
integrable  quadratic  differentials  on  $X$ by the pairing $\langle
\mu,  q\rangle =\int_{X} q\mu$. So the cotangent bundle $ T^*\mathcal
C  $  of $\mathcal C$ can be viewed as a suborbifold of $\mathcal Q_g
$,  the  moduli  space  of quadratic differentials. We will denote the total space of the cotangent bundle to the curve $\mathcal C$ by $\tilde{\mathcal C}$, and points of $\tilde{\mathcal C}$ by $(X, q)$, where $X$ is
a  Riemann  surface  and  $q$ is a quadratic differential on $X$. The
pullback  of  $\mathbb  L$  to  $\tilde{\mathcal C}$ will also be denoted by
$\mathbb L$. In the following, we will identify the cotangent bundle $\tilde{\mathcal C}$ with the tangent bundle by duality, so a quadratic differential will be seen as a tangent vector to $\mathcal C$.

We start by comparison of the  two  natural  metrics  on $\mathcal C$: the canonical hyperbolic
metric given by Riemann's uniformization theorem, and the induced Teichm\"uller metric. Both of these metrics are,
infinitesimally,  Finsler  metrics,  so  they  define  norms  on each
tangent space $T_X\mathcal C$ of $\mathcal C$.

\begin{lem}\label{comp}
Globally,  on  $\mathcal C$, the hyperbolic metric is larger than the
induced  \mbox{Teichm\"uller} metric, that is, the hyperbolic distance between
any two points is larger than the \mbox{Teichm\"uller} distance.

Infinitesimally, on each tangent space $T_X\mathcal C$, the unit ball
for  the  norm associated to the hyperbolic metric is included in the
unit ball for the norm associated to the induced Teichm\"uller metric.
\end{lem}

\begin{proof}
The  proof  is based on the notion of Kobayashi metric (cf \cite{H}).
The  canonical hyperbolic metric on $\mathcal C$ is by definition the
Kobayashi  metric  on  $\mathcal  C$,  and  by  Royden's theorem, the
Teichm\"uller  metric is the Kobayashi metric on $\mathcal M_g$. So the
statement  of  the  lemma results from the property of contraction of
the  (global  or  infinitesimal)  Kobayashi  metric for the inclusion
$\mathcal C \hookrightarrow \mathcal M_g$.
\end{proof}

Now we apply this lemma to the cotangent bundle $\tilde{\mathcal C}$, using the identification $\tilde{\mathcal C}\simeq T\mathcal C$.

\begin{cor}\label{in1}
Let  $\gamma(t)$  be  a  geodesic  on $\mathcal C$ for the hyperbolic
metric.   Let   us   denote by  $\gamma(\tau)$  the  same  curve
parameterized by the arc length for the Teichm\"uller metric restricted
to $\mathcal C$. The corresponding derivatives will be denoted by $\gamma'=\frac{\partial \gamma}{\partial t}$ and $\dot{\gamma}=\frac{\partial \gamma}{\partial \tau}$. Let $v$ be an element of $\mathbb L$ at $(\gamma(0),
\gamma'(0))=(X,  q/\Vert q \Vert_\textrm{hyp})\in \tilde {\mathcal C}
\subset  \mathcal  Q_g$.  Then the Lie derivatives of the norm of $v$
along $\gamma$ satisfy
\begin{equation}\label{met}
\left|\mathcal L_{\gamma'(0)}\log\Vert v\Vert\right|\leq
\left|\mathcal L_{\dot{\gamma}(0)}\log\Vert
v\Vert\right|
\end{equation}
\end{cor}

\begin{proof}
Note  that $\gamma'(0)$ and $\dot{\gamma}(0)$ are tangent vectors to
the  same  curve parametrized in two ways, at the same point, so they
are colinear:
$$
\gamma'(0)=\alpha \dot{\gamma}(0)\,.
$$
Since   $\gamma'(0)$  is  unitary  for  the  hyperbolic  metric,  and
$\dot{\gamma}(0)$  unitary  for  the  \mbox{Teichm\"uller}  metric,  by Lemma
\ref{comp},  $|\alpha|  \leq  1$. The conclusion follows by the chain
rule.
\end{proof}

Note  that Corollary~\ref{in1} is valid for any choice of the norm in
the  Hodge bundle provided the norm varies smoothly with respect to a
variation  of  a point in the base of the bundle. In the next section
we pass to a very special \textit{Hodge norm}.

\section{Variation of the Hodge norm}

The  natural  Hermitian  form  is  positive-definite  on  $H^{1,0}(X,
\mathbb   C)$,   so   it   induces   a   norm:
$$
\Vert  h\Vert^2_{1,0}=\frac{i}{2}\int_X  h\wedge\overline  h\,.
$$
Similarly, the intersection form is negative-definite on $H^{0,1}$ so
its  opposite  defines a norm $\Vert.\Vert_{0,1}$ on $ H^{0,1}$. Note
that  for  every  $h\in\mathbb  L^{1,0}$,  we  have  $\Vert \overline
h\Vert_  {0,1}=\Vert h\Vert_{1,0}$. We define the \textit{Hodge norm}
on  $H^1(X,  \mathbb  C)$  by  $\Vert v\Vert=\Vert h\Vert_{1,0}+\Vert
a\Vert_{0,1}$,  where  $h$ is the holomorphic part of $v$ and $a$ the
anti-holomorphic  part.  This  is  the  norm that we will consider on
$\mathbb L=\mathbb L^{1,0}\oplus\mathbb L^{0,1}$, by restriction.
From now on we consider only the Hodge norm.

The  second lemma gives a uniform bound for the variation of the Hodge norm in the direction of the Teichm\"uller flow  (for the definition of the Teichm\"uller flow, see e.g., \cite[Section 1]{Fo}).

\begin{lem}[G.  Forni]\label{in2}
Let  $v$ be a non trivial element of the fiber $H^1(X, \mathbb C)$ at
$(X,  q)\in  \tilde{\mathcal C}$. Then the Lie derivative of the Hodge norm  of  $v$  along  the  Teichm\"uller  flow  satisfies the following
inequality :
\begin{equation}\label{lie}
\left|\mathcal L \log\Vert v \Vert \right| \leq 1\,.
\end{equation}
Moreover, it  is  an  equality if and only if $q=\omega^ 2$ with $\omega\in H^{1,0}(\mathbb C)$ and
$$v\in
\textrm{Span}_{\mathbb C}(\omega)\oplus
\textrm{Span}_{\mathbb C}(\bar\omega)-\{0\}
$$
where $\textrm{Span}_{\mathbb C}(\omega)$ is the tautological bundle.
\end{lem}

The statement is the extension to the complex case of a lemma of G. Forni (see \cite[Lemma 2.1']{Fo}), reformulated in \cite[Cor. 2.1]{FMZ2}. The original statement holds in $H^1(\mathbb R)$ (with the Hodge norm), and by the Hodge representation theorem, it also holds in $H^{1,0}$, endowed with the norm $\Vert . \Vert _{1,0}$. By conjugation, we obtain the result in $H^{0,1}$ with the norm $\Vert.\Vert_ {0,1}$. Finally, note that inequality (\ref{lie}) is equivalent to the following: 
$$\left| \mathcal L \Vert v\Vert\right|\leq \Vert v\Vert. $$
So, applying this majoration to each component (holomorphic and antiholomorphic) of an element $v$ of $H^1(\mathbb C)$, endowed with the chosen norm $\Vert . \Vert$, we obtain the following inequalities: 
$$\left|\mathcal L \Vert v\Vert \right|= \left|\mathcal L\Vert h\Vert_{1,0}+\mathcal L \Vert a\Vert_{0,1}\right| \leq \left|\mathcal L\Vert h\Vert_{1,0}\right|+\left|\mathcal L \Vert a\Vert_{0,1}\right| \leq \Vert h\Vert_{1,0}+\Vert a\Vert_{0,1} = \Vert v\Vert,$$ so the result holds in $H^{1}(\mathbb C)$.

Note that there is another proof of inequality (\ref{lie}) in \cite[Lemma 6.10]{moller2}, in terms of curvature of the metrics. 

With this two lemmas we can achieve the proof of the theorem.

\section{Criterion in terms of Lyapunov exponents}

In this section we give an alternative version of the criterion, in terms of \mbox{Lyapunov} exponents. This version does not require any assumption on the signature of the pseudo-Hermitian form on the flat bundle $\mathbb L$, so it is more general, but it has less interest in practice, because Lyapunov exponents are harder to compute than orbifold degree. 

\begin{prop} 
Let $\mathcal C$ be a curve in the moduli space $\mathcal M_g$, with $\chi(\mathcal C)<0$, endowed with a flat subbundle $\mathbb L$ of rank $r\geq 2$ of the complex Hodge bundle, equivariant for the Gauss-Manin connection. Consider the Lyapunov exponents associated to the parallel transport of fibers of $\mathbb L$ along the geodesic flow given by the hyperbolic metric on $\mathcal C$. The absolute values of all these Lyapunov exponents are bounded above by $1$. If the bound is achieved, the curve $\mathcal C$ is a Teichm\"uller curve.
\end{prop}

\begin{proof}
Let us first explain where these Lyapunov exponents come from. 
Recall that $\mathcal C$ is endowed with a
canonical   hyperbolic   metric,  which  gives  us  a  geodesic  flow
$g_t^\textrm{hyp}$  on $T_1\mathcal C$, the unit tangent bundle of
$\mathcal  C$, and by duality, on the unit cotangent bundle $\tilde{\mathcal C}^{(1)}$. 
We  look  at the parallel transport of fibers of $\mathbb L$, endowed
with the Gauss-Manin connection, along this geodesic flow.

Let  $\nu$  be  the  Liouville measure  on $\tilde{\mathcal C}^{(1)}$. 

Since $\mathbb L$ inherits a variation of the Hodge structure from the
Hodge  bundle,  it has quasi-unipotent monodromy around any cusp (see
\cite[Th.  6.1]{schmid}).  So  there exists a finite unramified cover
$\hat{\mathcal  C}$  of  $\mathcal  C$,  such  that  the  pullback of
$\mathbb  L$ on $\hat{\mathcal C}$ has unipotent monodromy around any
cusp  of  $\hat{\mathcal  C}$. Passing to this finite cover preserves
the  ergodicity  of  the  geodesic  flow  that  we  consider on
$\hat{\mathcal  C}$ (because of the hyperbolic features of the geodesic flow and {\it Hopf's argument}, see e.g., Wilkinson's article \cite{wilkinson}),  and does not change the Lyapunov exponents. Actually, all  results  we  will obtain on $\mathcal C$ lift to $\hat{\mathcal C}$,  so  up  to passing to this cover, we will assume in the rest of this  paper that the monodromy of $\mathbb L$ is unipotent around any cusp of $\mathcal C$.

So with this assumption and thanks to the majoration given by Corollary \ref{in1} and Lemma \ref{in2}, the cocyle associated to the geodesic  flow  is log-integrable. The Oseledets theorem (see \cite{Oseledets}) can be applied to this complex cocycle. We denote by $\lambda_i$ the Lyapunov exponents and $E_{\lambda_i}$ the corresponding subspaces.

By definition every  vector  $v$  in  $E_{\lambda_i}(q)$  expands  with  the rate :
$$\lambda_i=\lim\limits_{t \to \infty}\frac{1}{t}\log
\|v(\gamma_q(t))\|,$$ where $\gamma_q$ is the geodesic for the hyperbolic metric starting at point $X$ in the direction $q$.

We can write :
\begin{equation}\label{inté}
\lambda_i=\lim\limits_{T \to \infty}
\frac{1}{T}\int_0^T \frac{d}{dt}\log\Vert v(\gamma_q(t))\Vert dt\,.
\end{equation}

So we have the following majoration :
\begin{equation}\label{in3} \lambda_i 
\leq   \lim\limits_{T  \to  \infty}  \frac{1}{T}\int_0^T \max_{v\in
\mathbb     L}\left(\frac{d}{dt}\log\Vert     v(\gamma_q(t))\Vert
\right) dt\,.
\end{equation}

By Birkhoff's theorem, we have : \begin{equation}\label{eq}\lim\limits_{T \to \infty} \frac{1}{T}\int_0^T \max_{v\in \mathbb L}\left(\frac{d}{dt}\log\Vert v(\gamma_q(t))\Vert \right) dt = \int_{\tilde{\mathcal C}^{(1)}} \max_{v\in \mathbb L}\left(\mathcal L _{\gamma_q'(0)}\log \Vert v(q)\Vert\right) d\nu(q).\end{equation}

It results from Corollary \ref{in1} and Lemma \ref{in2} that \begin{equation}\label{in4} \max_{v\in \mathbb L}\left\vert \mathcal L _{\gamma_q'(0)}\log \Vert v(q)\Vert\right\vert\leq 1,\end{equation} 
so \begin{equation}\label{in5} \int_{\tilde{\mathcal C}^{(1)}} \max_{v\in \mathbb L}\left( \mathcal L _{\gamma_q'(0)}\log \Vert v(q)\Vert\right) d\nu(q)\leq 1.\end{equation}

Similary, we have: \begin{equation}\label{in6}\lambda_i \geq \int_{\tilde{\mathcal C}^{(1)}} \min_{v\in \mathbb L}\left( \mathcal L _{\gamma_q'(0)}\log \Vert v(q)\Vert\right) d\nu(q)\geq -1.\end{equation}

So we obtain that $\vert\lambda_i\vert$ is bounded by $1$. Assume now that the bound is achieved for some index $i$. By symmetry of the spectrum (see \cite[Theorem 4]{FMZ3}), the exponent $-\lambda_i$ lies in the spectrum, so we can assume that $\lambda_i=1$. It means that inequalities (\ref{in3}) and (\ref{in5}) are in fact equalities. Since the measure $\nu$ is normalized, and the modulus of the integrand in (\ref{in5}) is $1$ at most (cf. (\ref{in4})), it is almost everywhere equal to $1$ and hence, by continuity, everywhere equal to $1$. It means that inequalities of Corollary \ref{in1} and Lemma \ref{in2} are equalities.

The first equality case in Corollary \ref{in1} implies that the two metrics (hyperbolic and Teichm\"uller) coincide on $\mathcal C$. Hence $\mathcal C$ is invariant by the Teichm\"uller flow, so it is a Teichm\"uller curve.

Let us denote by $v_1$ the element of $\mathbb L$ at point $(X,q)$ which maximizes the quantity $\mathcal L \log\Vert v\Vert\in[-1, 1]$. Similary, considering inequation (\ref{in6}) with $\lambda_i=-1$ gives an element $v_2$ which minimizes the same quantity. Clearly, $v_1$ and $v_2$ are independant. By Lemma \ref{in2}, we have $q=\omega^2$ and $\textrm{Span}_{\mathbb C}(v_1, v_2)=\textrm{Span}_{\mathbb C}(\omega, \overline{\omega})$. So we obtain the additional information that $\mathbb L$ contains the complex tautological bundle. In particular, $\mathbb L^{1,0}$ corresponds to $\textrm{Span}_{\mathbb C}(\omega)$.
\end{proof}

\section{End of the proof}

We will finish the proof of the theorem using that, with the additional assumption on the signature, there is only one non-negative Lyapunov exponent, which can be written in terms of degree of the line bundle $\mathbb L ^{1,0}$. So we will be able to conlude with the previous proposition.

 Since the pseudo-Hermitian form has signature
$(1,r-1)$ on $\mathbb L$, there is at most one positive Lyapunov exponent, denoted by $\lambda_1$ (see \cite[Theorem 4]{FMZ3}).

As  the  monodromy  is  unipotent  around  cusps,  the  degree of the
extended line bundle $\overline{\mathbb L^{1,0}} $ is the integral on
$\mathcal  C$  of the curvature form $\alpha$ (first Chern class), cf
\cite[Prop.  3.4]{Pe}.  Then  one has:  $$\int_{\mathcal  C}  \alpha=  \deg  \overline  {\mathbb
L^{1,0}}.$$

We  use  now  the formula for the sum of the Lyapunov exponents of an
invariant subbundle with respect to a geodesic flow defined  by  the  hyperbolic  metric on the curve $\mathcal{C}$. This formula  was  outlined  by  M.  Kontsevich  in \cite{K}, developed by
G.~Forni  in \cite{Fo}. I. Bouw and M. M\"oller suggested in~\cite{BM}
an  algebro-geometric interpretation of the numerator as the orbifold
degree of the associate line bundle.
As  it  was  mentioned in \cite[Section 2.3]{courtEKZ}, the result holds for
any  abstract  geodesic  flow.  So  here we apply this result for the
geodesic  flow  given  by  the  hyperbolic metric on $\mathcal C$: \begin{equation}\label{exp}
\lambda_1=-\frac{2\int_{\mathcal C}\alpha}{\chi(\mathcal C)}=-\frac{2\deg \overline{\mathbb L^{1,0}}}{\chi(\mathcal C)}.\end{equation}

This equality together with the previous proposition prove the second part of the theorem. 
The last statement of the theorem underlines the fact that any \mbox{Teichm\"uller} curve corresponding to a stratum of Abelian differentials admits a tautological bundle, which Lyapunov exponent is equal to 1.

This achieves the proof of the theorem.

\section{Acknowledgements}

The author is grateful to Quentin Gendron for helpful comments,
to \mbox{Martin M\"oller} and Giovanni Forni for valuable remarks, to Anton Zorich for formulation
of the problem, and to Carlos Matheus for careful reading the manuscript.

\bibliographystyle{amsalpha}

\end{document}